\documentclass{amsart}
\usepackage{latexsym}
\usepackage{amsfonts}
\usepackage{amssymb}
\usepackage{graphicx}

\newtheorem{theorem}{Theorem}[section]
\newtheorem{corollary}[theorem]{Corollary}
\newtheorem{lemma}[theorem]{Lemma}

\newtheorem{fact}{Fact}[section]

\theoremstyle{remark}
\newtheorem{remark}{Remark}[section]

\numberwithin{equation}{section}

\newcommand{\ba}{\mathfrak{a}}
\newcommand{\ab}{\mathfrak{b}}

\newcommand{\g}{\mathfrak{g}}
\newcommand{\gl}{\mathfrak{gl}}
\newcommand{\ah}{\mathfrak{h}}

\newcommand{\az}{\mathfrak{z}}

\newcommand{\CC}{\mathbb{C}}

\newcommand{\RR}{\mathbb{R}}

\newcommand{\ad}{\mathrm{ad}}
\newcommand{\Ad}{\mathrm{Ad}}

\newcommand{\Aut}{\mathrm{Aut}}

\newcommand{\im}{\mathrm{Im}}

\newcommand{\Rad}{\mathrm{Rad}}
\newcommand{\diff}{\frac{d}{dt}\Big|_{t=0}}
\newcommand{\GGG}{\overline{\overline{G}}}
\newcommand{\ggggg}{\overline{\overline{\g}}}
\newcommand{\GGGG}{\overline{\overline{G'}}}

\begin{document}

\title[Conjugacy classes in Lie groups]
{Topological properties of $\Ad$-semisimple conjugacy
classes in Lie groups}

\author{Jinpeng An}
\address{Department of Mathematics, ETH Zurich, 8092 Zurich, Switzerland}
\email{anjp@math.ethz.ch}

\keywords{Lie group, Lie algebra, conjugacy class, adjoint orbit.}

\subjclass[2000]{22E15; 17B05; 57S25.}

\begin{abstract}
We prove that $\Ad$-semisimple conjugacy classes in a connected
Lie group $G$ are closed embedded submanifolds of $G$. We also
prove that if $\alpha:H\rightarrow G$ is a homomorphism of
connected Lie groups such that the kernel of $\alpha$ is discrete
in $H$, then for an $\Ad$-semisimple conjugacy class $C$ in $G$,
every connected component of $\alpha^{-1}(C)$ is a conjugacy class
in $H$. Corresponding results for adjoint orbits in real Lie
algebras are also proved.
\end{abstract}

\maketitle


\section{Introduction}

Conjugacy classes in linear algebraic groups are extensively
studied and have had many applications (see, for example,
\cite{Hu}). We recall two well-known results concerning
topological properties of semisimple conjugacy classes in
algebraic groups below (see \cite[Proposition 18.2]{Hu0} as well
as the proofs of \cite[Theorem 3.8]{Hu} and \cite[Lemma
1.10]{PV}).

\begin{fact}\label{Borel-Harish-Chandra}
Let $\mathbf{G}$ be a linear algebraic group defined over an
algebraically closed field. Then semisimple conjugacy classes in
$\mathbf{G}$ are Zariski closed.
\end{fact}

\begin{fact}\label{Richardson}
Let $\mathbf{G}$ be a linear algebraic group defined over an
algebraically closed field, $\mathbf{H}$ a Zariski closed subgroup
of $\mathbf{G}$. Then for every semisimple conjugacy class $C$ in
$\mathbf{G}$, $C\cap \mathbf{H}$ is a finite union of conjugacy
classes in $\mathbf{H}$.
\end{fact}

As for Lie groups which are isomorphic to the groups of
$\RR$-points of certain algebraic groups defined over $\RR$, we
have the following result (see \cite[Proposition 10.1]{BH}).

\begin{fact}\label{BHC}
Let $\mathbf{G}\subset GL_n(\CC)$ be a linear algebraic group
defined over $\RR$, $G=\mathbf{G}\cap GL_n(\RR)$. Then semisimple
conjugacy classes in $G$ are closed with respect to the Hausdorff
topology.
\end{fact}

But for general Lie groups, there are relatively less known
results concerning properties of conjugacy classes, partly because
of the lack of algebraic instruments and the nonlinearity involved
in the problems.

The purpose of this paper is to prove the Lie-theoretic
counterparts of Facts \ref{Borel-Harish-Chandra} and
\ref{Richardson}. The corresponding results for adjoint orbits in
real Lie algebras are also obtained. Since the notion of
semisimplicity makes no sense for elements in general Lie groups
or general real Lie algebras, we introduce the following notions.

Let $G$ be a Lie group with Lie algebra $\g$. An element $g$ of
$G$ is \emph{$\Ad$-semisimple} if $\Ad(g)$ is semisimple in
$GL(\g)$. An element $X$ of $\g$ is \emph{$\ad$-semisimple} if
$\ad(X)$ is semisimple in $\gl(\g)$. A conjugacy class $C$ in $G$
is \emph{$\Ad$-semisimple} if a (hence every) element of $C$ is
$\Ad$-semisimple. An adjoint orbit $O$ in $\g$ is
\emph{$\ad$-semisimple} if a (hence every) element of $O$ is
$\ad$-semisimple.

The main results in this paper are the following two assertions.

\begin{theorem}\label{T:main1}
Let $G$ be a connected Lie group with Lie algebra $\g$. Then\\
(i) $\Ad$-semisimple conjugacy classes in $G$ are closed embedded
submanifolds of $G$;\\
(ii) $\ad$-semisimple adjoint orbits in $\g$ are closed embedded
submanifolds of $\g$.
\end{theorem}

\begin{theorem}\label{T:main2}
Let $\alpha:H\rightarrow G$ be a homomorphism of connected Lie
groups. Suppose $\ker(\alpha)$ is discrete in $H$. Let $C$ be an
$\Ad$-semisimple conjugacy class in $G$, and let $O$ be an
$\ad$-semisimple adjoint orbit in the Lie algebra $\g$ of $G$. Then\\
(i) Every connected component of $\alpha^{-1}(C)$ is a conjugacy class in $H$;\\
(ii) Every connected component of $(d\alpha)^{-1}(O)$ is an
adjoint orbit in the Lie algebra $\ah$ of $H$.
\end{theorem}

The following corollary of Theorem \ref{T:main2} is obvious.

\begin{corollary}\label{C:main2}
Let $G$ be a connected Lie group with Lie algebra $\g$, and let
$H$ be a connected closed subgroup of $G$ with Lie algebra $\ah$.
Then\\
(i) For every $\Ad$-semisimple conjugacy class $C$ in $G$, every
connected component of $C\cap H$ is a conjugacy class in $H$;\\
(ii) For every $\ad$-semisimple adjoint orbit $O$ in $\g$, every
connected component of $O\cap\ah$ is an adjoint orbit in $\ah$.
\end{corollary}

Theorem \ref{T:main1} generalizes Fact \ref{BHC}. The closedness
of conjugacy classes of finite order in connected Lie groups,
which is a special case of item (i) of Theorem \ref{T:main1}, was
proved in \cite{AW} using real analytic geometry.

The proofs of Facts \ref{Borel-Harish-Chandra} and
\ref{Richardson} are based on computations of dimensions of
subvarieties of $G$ and $\g$. But such computations do not help to
prove Theorems \ref{T:main1} and \ref{T:main2} because of the
facts: (1) the dimension of the boundary of an immersed
submanifold of a manifold is not necessarily strictly smaller than
the dimension of the submanifold; (2) the intersection of two
submanifolds of a manifold is not necessarily a submanifold.

Basically, our proofs of Theorems \ref{T:main1} and \ref{T:main2}
are Lie-theoretic. But the utilization of some facts of real
algebraic groups is essential, including the multiplicative Jordan
decomposition on real algebraic groups and Whitney's Theorem on
the finiteness of the number of connected components of a real
algebraic variety (see the proofs of Lemmas \ref{L:adss},
\ref{L:linear}, and \ref{L:ss} below).

The starting point of our proofs of Theorems \ref{T:main1} and
\ref{T:main2} is the following assertion.

\begin{theorem}\label{T:main3}
Let $G$ be a connected Lie group with Lie algebra $\g$,
$\rho:G\rightarrow GL(\mathcal{V})$ a representation of $G$ in a
finite dimensional real vector space $\mathcal{V}$. Suppose
$\ker(\rho)$ is discrete in $G$. Let $f\in\RR[\lambda]$ be a real
polynomial without
multiple roots in $\CC$. Then\\
(i) Every connected component of the subset
$$Z(\rho,f)=\{g\in G|f(\rho(g))=0\}$$ of $G$ is a closed embedded
submanifold of $G$, and is a conjugacy class in $G$;\\
(ii) Every connected component of the subset
$$\az(d\rho,f)=\{X\in\g|f(d\rho(X))=0\}$$ of $\g$ is a closed
embedded submanifold of $\g$, and is an adjoint orbit in $\g$.
\end{theorem}

For a conjugacy class $C$ in a connected Lie group $G$, and an
adjoint orbit $O$ in the Lie algebra $\g$ of $G$, define the sets
$$\Gamma(C)=g^{-1}C\cap Z(G), \quad g\in C,$$
$$\gamma(O)=(-X+O)\cap Z(\g), \quad X\in O,$$
where $Z(G)$ and $Z(\g)$ are the centers of $G$ and $\g$,
respectively. It can be shown that $\Gamma(C)$ and $\gamma(O)$ are
independent of the choices of $g\in C$ and $X\in O$, and that
$\Gamma(C)$ is a subgroup of $Z(G)$, $\gamma(O)$ is an additive
subgroup of $Z(\g)$.

The proofs of Theorems \ref{T:main1} and \ref{T:main2} are also
based on the following result.

\begin{theorem}\label{T:main4}
Let $G$ be a connected Lie group with Lie algebra $\g$. We have\\
(i) If $C$ is an $\Ad$-semisimple conjugacy class in $G$, then
$\Gamma(C)$ is a finite subgroup of $Z(G)$;\\
(ii) If $O$ is an $\ad$-semisimple adjoint orbit in $\g$, then
$\gamma(O)$ is trivial.
\end{theorem}

%

Now we give a sketch of the contents of the following sections. In
Section 2 we prove Theorem \ref{T:main3}. The key point is the
construction of a map $\varphi_g:\g=\ba_1\oplus\ba_2\rightarrow G$
for $g\in Z(\rho,f)$, which maps some neighborhood $U$ of $0\in\g$
onto a neighborhood of $g$ diffeomorphically, where
$\ba_1=\ker(1-\Ad(g))$, $\ba_2=\im(1-\Ad(g))$. Then one can show
that for $Y\in U$, $\varphi_g(Y)\in Z(\rho,f)$ if and only if
$Y\in\ba_2$ if and only if $\varphi_g(Y)$ is conjugate to $g$.
This idea is essentially due to N.-H. Neeb. The proof in the Lie
algebra level is similar.

After showing some basic properties of $\Gamma(C)$ and
$\gamma(O)$, Theorem \ref{T:main4} is proved in Section 3. We
first prove two special cases of item (i) of the theorem for which
$G$ is linear or semisimple using properties of real algebraic
groups, and then reduce the general case of item (i) to these two
cases using Levi's decomposition. The finiteness of $\Gamma(C)$ is
mainly deduced from Whitney's Theorem which states that the number
of connected components of a real algebraic variety is finite. The
proof of item (ii) makes use of item (i). One can also give a
direct proof of item (ii), which is similar to that of item (i).

Theorem \ref{T:main1} is proved in Section 4. We first use Theorem
\ref{T:main3} to prove that the projection of an $\Ad$-semisimple
conjugacy class $C$ in $G$ on $G/Z(G)_0$ is closed. This implies
that $C\cdot Z(G)_0$ is closed. By Theorem \ref{T:main4}, one can
show that $C\cdot Z(G)_0$ is a fiber bundle over
$Z(G)_0/(\Gamma(C)\cap Z(G)_0)$, and the fiber above the identity
is just $C$. So $C$ is closed. The proof of the closedness of
$\ad$-semisimple adjoint orbits is similar.

In Section 5 we prove Theorem \ref{T:main2}. If $G$ is linear and
$\alpha$ is injective, item (i) of Theorem \ref{T:main2} can be
easily deduced from Theorem \ref{T:main3}. We first prove the
theorem under the assumptions that $\Gamma(C)$ is trivial and
$\alpha$ is injective using this observation by considering the
adjoint group of $G$, and then deduce the general case from this.
The proof of item (ii) is similar but much easier.

The author is grateful to K.-H. Neeb and J.-K. Yu for valuable
conversations.


\section{Characterizations of conjugacy classes by polynomials}

In this section we prove the following theorem, which is the
staring point of the proofs of Theorems \ref{T:main1} and
\ref{T:main2}.

\begin{theorem}\label{T:polynomial}
Let $G$ be a connected Lie group with Lie algebra $\g$,
$\rho:G\rightarrow GL(\mathcal{V})$ a representation of $G$ in a
finite dimensional real vector space $\mathcal{V}$. Suppose
$\ker(\rho)$ is discrete. Let $f\in\RR[\lambda]$ be a real
polynomial without
multiple roots in $\CC$. Then\\
(i) Every connected component of the subset
$$Z(\rho,f)=\{g\in G|f(\rho(g))=0\}$$ of $G$ is a closed embedded
submanifold of $G$, and is a conjugacy class in $G$;\\
(ii) Every connected component of the subset
$$\az(d\rho,f)=\{X\in\g|f(d\rho(X))=0\}$$ of $\g$ is a closed
embedded submanifold of $\g$, and is an adjoint orbit in $\g$.
\end{theorem}

\begin{proof}
(i) Denote $Z=Z(\rho,f)$. Firstly, we note that $Z$ is invariant
under the conjugation of $G$. So $Z$ is the union of some
conjugacy classes in $G$.

Let $g\in Z$. By the definition of the set $Z$, $f(\rho(g))=0$.
Since $f$ has no multiple roots, $\rho(g)$ is semisimple. We claim
that $g$ is $\Ad$-semisimple. Indeed, since $\ker(\rho)$ is
discrete, the differential $d\rho:\g\rightarrow\gl(\mathcal{V})$
of $\rho$ is injective. So the action of $\Ad(g)$ on $\g$ is
equivalent to the action of $\Ad(\rho(g))|_{d\rho(\g)}$ on
$d\rho(\g)$. Since $\rho(g)$ acts semisimply on $\mathcal{V}$,
$\Ad(\rho(g))$ acts semisimply on $\gl(\mathcal{V})$, and then
$\Ad(\rho(g))|_{d\rho(\g)}$ acts semisimply on $d\rho(\g)$. This
verifies the claim.

Denote $\ba_1=\ker(1-\Ad(g))$, $\ba_2=\im(1-\Ad(g))$. Since
$\Ad(g)$ is semisimple, $\g=\ba_1\oplus\ba_2$. Define a map
$\varphi_g:\ba_1\oplus\ba_2\rightarrow G$ by
$$\varphi_g(Y_1,Y_2)=e^{Y_2}e^{Y_1}ge^{-Y_2}, \quad \quad Y_1\in\ba_1, Y_2\in\ba_2.$$
Then it is easy to compute the differential
$(d\varphi_g)_{(0,0)}:\ba_1\oplus\ba_2\rightarrow T_gG$ of
$\varphi_g$ at $(0,0)$ as
$$(d\varphi_g)_{(0,0)}(Y_1,Y_2)=(dr_g)_e(Y_1+(1-\Ad(g))Y_2),$$ where $r_g$
is the right translation on $G$ induced by $g$. Since $\Ad(g)$ is
semisimple, the restriction of $1-\Ad(g)$ on $\ba_2=\im(1-\Ad(g))$
is a linear automorphism. Hence $(d\varphi_g)_{(0,0)}$ is a linear
isomorphism. By the Implicit Function Theorem, there exist an open
neighborhood $U_1$ of $0\in\ba_1$ and an open neighborhood $U_2$
of $0\in\ba_2$ such that the restriction of $\varphi_g$ to
$U_1\times U_2\subset\ba_1\oplus\ba_2$ is a diffeomorphism onto an
open neighborhood $U=\varphi_g(U_1\times U_2)$ of $g\in G$.

Define a map $\alpha_g:\ba_1\rightarrow\gl(\mathcal{V})$ by
$$\alpha_g(Y_1)=f(\rho(e^{Y_1}g)).$$ We claim that $\alpha_g$ is an
immersion at $0\in\ba_1$. Indeed, we have
\begin{align*}
&(d\alpha_g)_0(Y_1)\\
=&\diff\alpha_g(tY_1)\\
=&\diff f(e^{td\rho(Y_1)}\rho(g))\\
=&d\rho(Y_1)\rho(g)f'(\rho(g)),
\end{align*}
where $f'$ is the derivative of $f$. Here the last step holds
because $e^{td\rho(Y_1)}$ commutes with $\rho(g)$. Since $f$ has
no multiple roots, $(f,f')=1$. So there exist polynomials $r,s$
such that $fr+f's=1$. Substitute $\rho(g)$ for the indeterminate
in this equality and notice that $f(\rho(g))=0$, we get
$f'(\rho(g))s(\rho(g))=1$. So $f'(\rho(g))$ is invertible. Since
$\rho(g)$ is also invertible and $d\rho$ is injective,
$(d\alpha_g)_0(Y_1)=0$ implies $Y_1=0$. Hence $\alpha_g$ is an
immersion at $0\in\ba_1$. Thus, shrinking $U_1$ if necessary, we
may assume that $\alpha_g|_{U_1}$ is injective.

Now for $Y_1\in U_1$, $Y_2\in U_2$, we have
\begin{align*}
&f(\rho(\varphi_g(Y_1,Y_2)))\\
=&f(\rho(e^{Y_2})\rho(e^{Y_1}g)\rho(e^{Y_2})^{-1})\\
=&\rho(e^{Y_2})f(\rho(e^{Y_1}g))\rho(e^{Y_2})^{-1}\\
=&\rho(e^{Y_2})\alpha_g(Y_1)\rho(e^{Y_2})^{-1}.
\end{align*}
So $\varphi_g(Y_1,Y_2)\in Z\Leftrightarrow Y_1=0$, that is,
$$Z\cap U=\varphi_g(\{0\}\times U_2)=\{e^{Y_2}ge^{-Y_2}|Y_2\in U_2\}.$$
This shows that every connected component of $Z$ is an embedded
submanifold of $G$, which is necessarily closed by the definition
of $Z$, and that every conjugacy class contained in $Z$ is open in
$Z$. But the connectedness of $G$ implies that conjugacy classes
are connected. Hence every conjugacy class contained in $Z$ is in
fact a connected component of $Z$. This proves (i).

(ii) Similar to the proof of (i), the set $\az=\az(d\rho,f)$ is
the union of some adjoint orbits in $\g$. Let $X\in\az$. Then
$d\rho(X)$ and $\ad(X)$ are semisimple. Denote
$\ab_1=\ker(\ad(X))$, $\ab_2=\im(\ad(X))$. Then
$\g=\ab_1\oplus\ab_2$. Define a map
$\psi_X:\ab_1\oplus\ab_2\rightarrow\g$ by
$$\psi_X(W_1,W_2)=\Ad(e^{W_2})(X+W_1), \quad \quad W_1\in\ab_1, W_2\in\ab_2.$$ Then
$$(d\psi_X)_{(0,0)}(W_1,W_2)=W_1-\ad(X)(W_2).$$  Hence $(d\psi_X)_{(0,0)}$ is a linear
isomorphism, and then there exist an open neighborhood $V_1$ of
$0\in\ab_1$ and an open neighborhood $V_2$ of $0\in\ab_2$ such
that the restriction of $\psi_X$ to $V_1\times
V_2\subset\ab_1\oplus\ab_2$ is a diffeomorphism onto an open
neighborhood $V=\psi_X(V_1\times V_2)$ of $X\in\g$.

Define $\beta_X:\ab_1\rightarrow\gl(\mathcal{V})$ by
$$\beta_X(W_1)=f(d\rho(X+W_1)).$$ Then
\begin{align*}
&(d\beta_X)_0(W_1)\\
=&\diff\beta_X(tW_1)\\
=&\diff f(d\rho(X)+td\rho(W_1))\\
=&d\rho(W_1)f'(d\rho(X)).
\end{align*}
Similar to the proof of (i), we can prove $f'(d\rho(X))$ is
invertible. So $(d\beta_X)_0(W_1)=0$ implies $W_1=0$, that is,
$\beta_X$ is an immersion at $0\in\ab_1$. Shrinking $V_1$ if
necessary, we may assume that $\beta_X|_{V_1}$ is injective.

Now for $W_1\in V_1$, $W_2\in V_2$, we have
\begin{align*}
&f(d\rho(\psi_X(W_1,W_2)))\\
=&f(d\rho(\Ad(e^{W_2})(X+W_1)))\\
=&f(\rho(e^{W_2})d\rho(X+W_1)\rho(e^{W_2})^{-1})\\
=&\rho(e^{W_2})f(d\rho(X+W_1))\rho(e^{W_2})^{-1}\\
=&\rho(e^{W_2})\beta_X(W_1)\rho(e^{W_2})^{-1}.
\end{align*}
So $\psi_X(W_1,W_2)\in\az\Leftrightarrow W_1=0$, that is,
$$\az\cap V=\psi_X(\{0\}\times V_2)=\{\Ad(e^{W_2})(X)|W_2\in V_2\}.$$
Then an argument similar to the proof of (i) shows that every
connected component of $\az$ is a closed embedded submanifold of
$\g$, and is an adjoint orbit. This proves (ii).
\end{proof}

\begin{corollary}\label{C:polynomial}
Let $G$ be a connected Lie group with Lie algebra $\g$, and let
$\rho:G\rightarrow GL(\mathcal{V})$ be a representation of $G$ in
a finite dimensional real vector space $\mathcal{V}$. Suppose
$\ker(\rho)$ is discrete. We have\\
(i) If $C$ is a conjugacy class in $G$ such that $\rho(C)$
contains a semisimple element $A$ of $GL(\mathcal{V})$, then $C$
is a closed embedded submanifold of $G$, and is a connected
component of the set
$$Z=\{g\in G|f(\rho(g))=0\},$$ where $f$ is the minimal polynomial
of $A$;\\
(ii) If $O$ is an adjoint orbit in $\g$ such that $d\rho(O)$
contains a semisimple element $B$ of $\gl(\mathcal{V})$, then $O$
is a closed embedded submanifold of $\g$, and is a connected
component of the set
$$\az=\{X\in\g|p(d\rho(X))=0\},$$ where $p$ is the minimal polynomial
of $B$.\qed
\end{corollary}


\section{Finiteness of $\Gamma(C)$
}

Let $G$ be a connected Lie group with Lie algebra $\g$. Let $C$ be
a conjugacy class in $G$, and let $O$ be an adjoint orbit in $\g$.
The subset $\Gamma(C)=\Gamma_G(C)$ of the center $Z(G)$ of $G$ is
defined by
$$\Gamma_G(C)=g^{-1}C\cap Z(G), \quad g\in C.$$ The
subset $\gamma(O)=\gamma_\g(O)$ of the center $Z(\g)$ of $\g$ is
defined by
$$\gamma_\g(O)=(-X+O)\cap Z(\g), \quad X\in O.$$ For convenience,
denote $$\Gamma_0(C)=\Gamma(C)\cap Z(G)_0=g^{-1}C\cap Z(G)_0,
\quad g\in C,$$ where $Z(G)_0$ is the identity component of
$Z(G)$. In this section we prove some properties of $\Gamma(C)$
and $\gamma(O)$, especially the finiteness of $\Gamma(C)$ and the
triviality of $\gamma(O)$ under the $\Ad$-semisimplicity or
$\ad$-semisimplicity condition.

\begin{lemma}\label{L:independent}
Let $G$ be a connected Lie group, $C$ a conjugacy class in $G$,
$O$ an adjoint orbit in the Lie algebra $\g$ of $G$. Then\\
(i) $\Gamma(C)$ is independent of the choice of the element $g\in
C$ defining it;\\
(ii) $\gamma(O)$ is independent of the choice of the element $X\in
O$ defining it.
\end{lemma}

\begin{proof}
(i) Let $g_1,g_2\in C$. Then $g_1=hg_2h^{-1}$ for some $h\in G$.
Hence we have
\begin{align*}
&g_1^{-1}C\cap Z(G)=hg_2^{-1}h^{-1}C\cap
Z(G)=hg_2^{-1}(h^{-1}Ch)h^{-1}\cap
Z(G)\\
=&hg_2^{-1}Ch^{-1}\cap Z(G)=h(g_2^{-1}C\cap
Z(G))h^{-1}=g_2^{-1}C\cap Z(G).
\end{align*}
This proves (i).

(ii) Let $X_1,X_2\in O$. Then $X_1=\Ad(g)X_2$ for some $g\in G$.
Hence
\begin{align*}
&(-X_1+O)\cap Z(\g)=(-\Ad(g)X_2+O)\cap Z(\g)\\
=&\Ad(g)(-X_2+O)\cap Z(\g)=\Ad(g)((-X_2+O)\cap
Z(\g))\\
=&(-X_2+O)\cap Z(\g).
\end{align*}
This proves (ii).
\end{proof}

For an element $g$ in a connected Lie group $G$, we denote by
$Z_G(g)$ the centralizer of $g$ in $G$, and denote
$$N_G(g)=\{h\in G|g^{-1}hgh^{-1}\in Z(G)\}.$$
$N_G(g)$ is a closed subgroup of $G$ containing $Z_G(g)$. In fact,
if we let $\pi:G\rightarrow G/Z(G)$ be the quotient homomorphism,
then $N_G(g)=\pi^{-1}(Z_{G/Z(G)}(\pi(g)))$. Similarly, for an
element $X$ in the Lie algebra $\g$ of $G$, denote by $Z_G(X)$ the
centralizer of $X$ in $G$, and denote
$$N_G(X)=\{h\in G|-X+\Ad(h)X\in Z(\g)\}.$$ Then
$N_G(X)=\pi^{-1}(Z_{G/Z(G)}(d\pi(X)))$ is a
closed subgroup of $G$ containing $Z_G(X)$.

\begin{lemma}\label{L:subgroup}
Let $G$ be a connected Lie group with Lie algebra $\g$, $C$ a
conjugacy class in $G$, $O$ an adjoint orbit in $\g$. Then\\
(i) $\Gamma(C)$ is a Lie subgroup of $Z(G)$, and is isomorphic to
$N_G(g)/Z_G(g)$ for every $g\in C$;\\
(ii) $\gamma(O)$ is a Lie subgroup of the vector group $Z(\g)$,
and is isomorphic to $N_G(X)/Z_G(X)$ for every $X\in O$.
\end{lemma}

\begin{proof}
(i) Let $g\in C$. Define a smooth map $\alpha:N_G(g)\rightarrow
Z(G)$ by
$$\alpha(h)=g^{-1}hgh^{-1}.$$ We claim
that $\alpha$ is a homomorphism of Lie groups. Indeed, let
$h_1,h_2\in N_G(g)$, then
\begin{align*}
&\alpha(h_1)\alpha(h_2)=(g^{-1}h_1gh_1^{-1})(g^{-1}h_2gh_2^{-1})\\
=&g^{-1}h_1g(g^{-1}h_2gh_2^{-1})h_1^{-1}=g^{-1}(h_1h_2)g(h_1h_2)^{-1}\\
=&\alpha(h_1h_2).
\end{align*}
It is obvious that the kernel of $\alpha$ is $Z_G(g)$, and the
image of $\alpha$ is $\Gamma(C)=g^{-1}C\cap Z(G)$. So $\Gamma(C)$
is a Lie subgroup of $Z(G)$, and is isomorphic to $N_G(g)/Z_G(g)$.

(ii) Let $X\in O$. Define $\beta:N_G(X)\rightarrow Z(\g)$ by
$$\beta(h)=-X+\Ad(h)X.$$ For $g_1,g_2\in N_G(X)$, we have
\begin{align*}
&\beta(g_1g_2)=-X+\Ad(g_1g_2)X\\
=&(-X+\Ad(g_1)X)+(-\Ad(g_1)X+\Ad(g_1)\Ad(g_2)X)\\
=&\beta(g_1)+\Ad(g_1)(-X+\Ad(g_2)X)\\
=&\beta(g_1)+\Ad(g_1)\beta(g_2)=\beta(g_1)+\beta(g_2).
\end{align*}
So $\beta$ is a homomorphism of Lie groups. The kernel of $\beta$
is $Z_G(X)$, the image of $\beta$ is $\gamma(O)=(-X+O)\cap Z(\g)$.
So $\gamma(O)$ is a Lie subgroup of $Z(\g)$, and is isomorphic to
$N_G(X)/Z_G(X)$.
\end{proof}

For a connected Lie group $G$ with Lie algebra $\g$ and an adjoint
orbit $O$ in $\g$, $\exp(O)$ is a conjugacy class in $G$.
$\gamma(O)$ and $\Gamma(\exp(O))$ have the following relation.

\begin{lemma}\label{L:relation}
Let $G$ be a connected Lie group with Lie algebra $\g$, $O$ an
adjoint orbit in $\g$. Then
$\exp(\gamma(O))\subset\Gamma_0(\exp(O))$.
\end{lemma}

\begin{proof}
Let $X\in O$. If $Y\in\gamma(O)$, then there exists $h\in G$ such
that $Y=-X+\Ad(h)X$. Since $Y\in Z(\g)$,
$he^Xh^{-1}=e^{\Ad(h)X}=e^{X+Y}=e^Xe^Y$. So
$e^Y=e^{-X}he^Xh^{-1}\in e^{-X}\exp(O)\cap
Z(G)_0=\Gamma_0(\exp(O))$. This shows
$\exp(\gamma(O))\subset\Gamma_0(\exp(O))$.
\end{proof}

Let $\pi:G\rightarrow G'$ be a covering homomorphism of Lie
groups. Then for a conjugacy class $C$ in $G$, $\pi(C)$ is a
conjugacy class in $G'$. The next lemma relates $\Gamma_G(C)$ with
$\Gamma_{G'}(\pi(C))$.

\begin{lemma}\label{L:invariant}
Let $\pi:G\rightarrow G'$ be a covering homomorphism of connected
Lie groups, $C$ a conjugacy class in $G$. Then
$\pi(\Gamma_G(C))=\Gamma_{G'}(\pi(C))$.
\end{lemma}

\begin{proof}
First we claim that $Z(G)=\pi^{-1}(Z(G'))$. Indeed, let $z\in
\pi^{-1}(Z(G'))$, and let $\alpha:G\rightarrow G$ be the map
defined by $\alpha(h)=hzh^{-1}z^{-1}$. Then
$\alpha(G)\subset\ker(\pi)$. Since $\alpha(G)$ is connected
containing the identity $e$ of $G$, and $\ker(\pi)$ is discrete,
we have $\alpha(G)=\{e\}$. So $z\in Z(G)$. This shows
$\pi^{-1}(Z(G'))\subset Z(G)$. It is obvious that
$Z(G)\subset\pi^{-1}(Z(G'))$. Hence $Z(G)=\pi^{-1}(Z(G'))$. Now we
choose a $g\in C$, then
\begin{align*}
&\pi(\Gamma_G(C))=\pi(g^{-1}C\cap Z(G))=\pi(g^{-1}C\cap\pi^{-1}(Z(G')))\\
=&\pi(g^{-1}C)\cap Z(G')=\pi(g)^{-1}\pi(C)\cap
Z(G')=\Gamma_{G'}(\pi(C)).
\end{align*}
\end{proof}

The following lemma demonstrates a rough understanding of
$\Gamma(C)$ and $\gamma(O)$ under the semisimplicity assumptions.

\begin{lemma}\label{L:rough}
Let $G$ be a connected Lie group with Lie algebra $\g$. We have\\
(i) If $C$ is an $\Ad$-semisimple conjugacy class in $G$, $g\in
C$, then the Lie algebras of $N_G(g)$ and $Z_G(g)$ coincide, and
$\Gamma(C)$ is a $0$-dimensional Lie subgroup of $Z(G)$;\\
(ii) If $O$ is an $\ad$-semisimple adjoint orbit in $\g$, $X\in
O$, then the Lie algebras of $N_G(X)$ and $Z_G(X)$ coincide, and
$\gamma(O)$ is a $0$-dimensional Lie subgroup of the vector group
$Z(\g)$.
\end{lemma}

\begin{proof}
(i) Since $Z_G(g)\subset N_G(g)$, to prove their Lie algebras
coincide, it is sufficient to show that for every $X$ in the Lie
algebra of $N_G(g)$, $X$ belongs to the Lie algebra of $Z_G(g)$.
For such an $X$, we have $g^{-1}e^{tX}ge^{-tX}\in Z(G)$ for every
$t\in\RR$. So
$e^{tX}e^{-t\Ad(g)X}=g(g^{-1}e^{tX}ge^{-tX})g^{-1}\in Z(G)$. This
implies that $(1-\Ad(g))X$ belongs to the Lie algebra of $Z(G)$,
and then $(1-\Ad(g))^2X=0$. Since $C$ is $\Ad$-semisimple,
$\Ad(g)$ is semisimple. So we in fact have $(1-\Ad(g))X=0$. But
the Lie algebra of $Z_G(g)$ is $\ker(1-\Ad(g))$. So $X$ belongs to
the Lie algebra of $Z_G(g)$. Hence the Lie algebras of $N_G(g)$
and $Z_G(g)$ coincide. As the image of the homomorphism $\alpha$
constructed in the proof of Lemma \ref{L:subgroup}, $\Gamma(C)$ is
a $0$-dimensional Lie subgroup of $Z(G)$.

(ii) Similar to the proof of (i), let $Y$ be an element of the Lie
algebra of $N_G(X)$. Then $-X+\Ad(e^{tY})X\in Z(\g)$ for every
$t\in\RR$. This implies that $\ad(Y)X\in Z(\g)$. So
$\ad(X)^2Y=-\ad(X)(\ad(Y)X)=0$. Since $X$ is $\ad$-semisimple,
$\ad(X)Y=0$. This shows that $Y$ belongs to the Lie algebra of
$Z_G(X)$. So the Lie algebras of $N_G(X)$ and $Z_G(X)$ coincide,
and $\gamma(O)$ is a $0$-dimensional Lie subgroup of $Z(\g)$.
\end{proof}

\begin{remark}
We only need the discreteness of $\Gamma(C)$ in $Z(G)$ in the
proofs of items (i) of Theorems \ref{T:main1} and \ref{T:main2}.
By Lemma \ref{L:rough}, $\Gamma(C)$ is $0$-dimensional when $C$ is
$\Ad$-semisimple. But this does not imply that $\Gamma(C)$ is
discrete in $Z(G)$.  To get the discreteness of $\Gamma(C)$, we
have to show that it is finite. If fact, if $\Gamma(C)$ could be
infinite for some connected Lie group $G$ and some
$\Ad$-semisimple conjugacy class $C$ in $G$, we would easily
construct a discrete central subgroup $D$ of $G\times\RR$ such
that $\Gamma_{(G\times\RR)/D}(\pi(C))$ is not discrete, where
$\pi:G\rightarrow(G\times\RR)/D$ is the covering homomorphism.
\end{remark}

The remaining of this section is devoted to the proof of Theorem
\ref{T:main4}. Some results on real algebraic groups are needed.
For convenience, we understand the Zariski topology on $GL_n(\RR)$
as the topology for which a closed set is the set of common zeros
of a family of real polynomial functions on $GL_n(\RR)$ with
indeterminates $g_{ij}$ $(1\leq i,j\leq n)$ and $\frac{1}{\det
g}$, where $g=(g_{ij})\in GL_n(\RR)$. It is obvious that if $G$ is
a Lie subgroup of $GL_n(\RR)$, then the Zariski closure $\GGG$ of
$G$ is also a Lie subgroup of $GL_n(\RR)$.

\begin{lemma}\label{L:adss}
Let $G$ be a connected Lie subgroup of $GL_n(\RR)$ for some $n$,
$\GGG$ the Zariski closure of $G$ in $GL_n(\RR)$. If $g\in G$ is
$\Ad$-semisimple in $G$, then it is $\Ad$-semisimple in $\GGG$.
\end{lemma}

\begin{proof}
Let $g=g_sg_u$ be the multiplicative Jordan decomposition of $g$
in $GL_n(\RR)$, where $g_s$ is semisimple, $g_u$ is unipotent. It
is well known that $g_s,g_u\in\GGG$ (see, for example,
\cite[Chapter 1, Section 4]{Bo}). Then
$\Ad_{\GGG}(g)=\Ad_{\GGG}(g_s)\cdot\Ad_{\GGG}(g_u)$ is the
multiplicative Jordan decomposition of $\Ad_{\GGG}(g)$ in
$GL(\ggggg)$, where $\ggggg$ is the Lie algebra of $\GGG$. Since
the Lie algebra $\g$ of $G$ in invariant under $\Ad_{\GGG}(g)$, it
is also invariant under $\Ad_{\GGG}(g_s)$ and $\Ad_{\GGG}(g_u)$.
So $\Ad(g)=\Ad_{\GGG}(g_s)|_\g\cdot\Ad_{\GGG}(g_u)|_\g$ is the
multiplicative Jordan decomposition of $\Ad(g)$ in $GL(\g)$. But
by the assumption, $\Ad(g)$ is semisimple. So
$\Ad_{\GGG}(g_u)|_\g=0$. This implies that $G\subset
Z_{\GGG}(g_u)$. Since $Z_{\GGG}(g_u)$ is Zariski closed, we have
$\GGG\subset Z_{\GGG}(g_u)$, that is, $g_u\in Z(\GGG)$. So
$\Ad_{\GGG}(g_u)=1$, and then $\Ad_{\GGG}(g)=\Ad_{\GGG}(g_s)$ is
semisimple, that is, $g$ is $\Ad$-semisimple in $\GGG$.
\end{proof}

\begin{lemma}\label{L:linear}
Let $G$ be a connected Lie subgroup of $GL_n(\RR)$ for some $n$,
$C$ an $\Ad$-semisimple conjugacy class in $G$. Then $\Gamma(C)$
is a finite subgroup of $Z(G)$.
\end{lemma}

\begin{proof}
Let $\GGG$ be the Zariski closure of $G$ in $GL_n(\RR)$, and let
$C'$ be the conjugacy class in $\GGG$ containing $C$. Choose a
$g\in C$. Since $g$ is $\Ad$-semisimple in $G$, by Lemma
\ref{L:adss}, $g$ is $\Ad$-semisimple in $\GGG$. by Lemma
\ref{L:rough}, the Lie algebras of $N_{\GGG}(g)$ and $Z_{\GGG}(g)$
coincide. Since $N_{\GGG}(g)$ can be expressed as
$$N_{\GGG}(g)=\{h\in\GGG|(g^{-1}hgh^{-1})x=x(g^{-1}hgh^{-1}), \forall
x\in\GGG\},$$ which is algebraic, by Whitney's Theorem \cite{Wh},
$N_{\GGG}(g)$ has finitely many connected components. So as a
quotient group of the component group of $N_{\GGG}(g)$,
$N_{\GGG}(g)/Z_{\GGG}(g)$ is finite. Hence $\Gamma_{\GGG}(C')\cong
N_{\GGG}(g)/Z_{\GGG}(g)$ is finite.

We claim that $Z(G)\subset Z(\GGG)$. Indeed, if $z\in Z(G)$, then
$Z_{GL_n(\RR)}(z)$ is an algebraic subgroup of $GL_n(\RR)$
containing $G$. So $Z_{GL_n(\RR)}(z)$ contains $\GGG$, that is,
$z\in Z(\GGG)$. This shows $Z(G)\subset Z(\GGG)$. Now we have
$$\Gamma_G(C)=g^{-1}C\cap Z(G)\subset g^{-1}C'\cap Z(\GGG)=\Gamma_{\GGG}(C').$$
Hence $\Gamma_G(C)$ is finite. This proves the lemma.
\end{proof}

\begin{lemma}\label{L:ss}
Let $G$ be a connected semisimple Lie group, $C$ an
$\Ad$-semisimple conjugacy class in $G$. Then $\Gamma(C)$ is a
finite subgroup of $Z(G)$.
\end{lemma}

\begin{proof}
Let $\Aut(\g)$ be the automorphism group of the Lie algebra $\g$
of $G$. Since $$\Aut(\g)=\{A\in GL(\g)|f_{X,Y}(A)=0, \forall
X,Y\in\g\},$$ where $$f_{X,Y}(A)=A[X,Y]-[AX,AY]$$ is algebraic,
$\Aut(\g)$ is an algebraic subgroup of $GL(\g)$. Choose $g\in C$.
Then $\Ad(g)\in\Aut(\g)$. By Whitney's Theorem,
$Z_{\Aut(\g)}(\Ad(g))$ has finitely many connected components.
Since $G$ is semisimple, $\Ad(G)$ is the identity component of
$\Aut(\g)$. So $Z_{\Ad(G)}(\Ad(g))=Z_{\Aut(\g)}(\Ad(g))\cap\Ad(G)$
has finitely many connected components. Since the kernel of the
epimorphism $\Ad:G\rightarrow\Ad(G)$ is $Z(G)$, which is discrete,
$N_G(g)/Z(G)\cong Z_{\Ad(G)}(\Ad(g))$ has finitely many connected
components.

On the other hand, we have
$$\Gamma(C)\cong N_G(g)/Z_G(g)\cong(N_G(g)/Z(G))/(Z_G(g)/Z(G)).$$
By Lemma \ref{L:rough}, the Lie algebras of $N_G(g)$ and $Z_G(g)$
coincide. So the Lie algebras of $N_G(g)/Z(G)$ and $Z_G(g)/Z(G)$
coincide. We have shown that $N_G(g)/Z(G)$ has finitely many
connected components. So
$\Gamma(C)\cong(N_G(g)/Z(G))/(Z_G(g)/Z(G))$ is finite.
\end{proof}

Now we can prove the finiteness of $\Gamma(C)$ and the triviality
of $\gamma(O)$ under the semisimplicity assumptions.

\begin{theorem}\label{T:finite}
Let $G$ be a connected Lie group with Lie algebra $\g$. We have\\
(i) If $C$ is an $\Ad$-semisimple conjugacy class in $G$, then
$\Gamma(C)$ is a finite subgroup of $Z(G)$;\\
(ii) If $O$ is an $\ad$-semisimple adjoint orbit in $\g$, then
$\gamma(O)$ is trivial.
\end{theorem}

\begin{proof}
(i) By Lemma \ref{L:invariant}, we may assume that $G$ is simply
connected. Let $R=\Rad(G)$. By Levi's Theorem, there is a
connected semisimple subgroup $L$ of $G$ such that $G=R\rtimes L$.
Note that $R$ and $L$ are simply connected.

We first prove that $R\cap\Gamma(C)$ is finite. Let $\Lambda(L)$
be the linearizer of $L$ (by definition, $\Lambda(L)$ is the
intersection of the kernels of all finite dimensional
representations of $L$). By considering the adjoint representation
of $L$ in the Lie algebra of $G$, we know that $\Lambda(L)\subset
Z(G)$. Since $L/\Lambda(L)$ admits a finite dimensional faithful
representation (see \cite[Chapter 5, Section 3, Theorem 8]{OV}),
by a theorem of Harish-Chandra \cite{HC}, $G/\Lambda(L)\cong
R\rtimes(L/\Lambda(L))$ admits a finite dimensional faithful
representation. Since $C$ is $\Ad$-semisimple in $G$, $\pi(C)$ is
$\Ad$-semisimple in $G/\Lambda(L)$, where $\pi:G\rightarrow
G/\Lambda(L)$ is the quotient homomorphism. By Lemma
\ref{L:linear}, $\Gamma_{G/\Lambda(L)}(\pi(C))$ is finite. Since
$\Lambda(L)$ is discrete, by Lemma \ref{L:invariant},
$\pi(\Gamma(C))=\Gamma_{G/\Lambda(L)}(\pi(C))$ is finite. Since
$R\cap\Lambda(L)$ is trivial, the restriction of $\pi$ to
$R\cap\Gamma(C)$ is injective. So $R\cap\Gamma(C)$ is finite.

Now consider the quotient homomorphism $\alpha:G\rightarrow G/R$.
Since $\alpha(C)$ is $\Ad$-semisimple in $G/R$, by Lemma
\ref{L:ss}, $\Gamma_{G/R}(\alpha(C))$ is finite. But the kernel of
the homomorphism $\alpha|_{\Gamma(C)}:\Gamma(C)\rightarrow
\Gamma_{G/R}(\alpha(C))$ is $R\cap\Gamma(C)$, which we have shown
is finite. So $\alpha(C)$ is finite. This proves (i).

(ii) We may assume that $G$ is simply connected. Since $O$ is
$\ad$-semisimple, $\exp(O)$ is an $\Ad$-semisimple conjugacy class
in $G$. By item (i) of the theorem, $\Gamma(\exp(O))$ is finite.
So $\Gamma_0(\exp(O))=\Gamma(\exp(O))\cap Z(G)_0$ is a finite
subgroup of $Z(G)_0$. But the simple connectedness of $G$ implies
that $Z(G)_0$ is simply connected (see \cite[Corollary
3.18.6]{Va}), which is isomorphic to a vector group. So
$\Gamma_0(\exp(O))$ is in fact trivial. By Lemma \ref{L:relation},
$\exp(\gamma(O))$ is trivial. But the simple connectedness of
$Z(G)_0$ implies that the restriction of the exponential map to
$Z(\g)$ is injective. In particular, $\exp|_{\gamma(O)}$ is
injective. So $\gamma(O)$ is trivial.
\end{proof}

\begin{remark}
Our proof of item (ii) of Theorem \ref{T:finite} is based on item
(i) of that theorem. But one can also give a direct proof of item
(ii). To do this, one can embed $\g$ into some $\gl_n(\RR)$ using
Ado's Theorem, consider the connected Lie subgroup $G'$ of
$GL_n(\RR)$ with Lie algebra $\g$, and then consider the Zariski
closure $\GGGG$ of $G'$. In this course one need a result similar
to Lemma \ref{L:adss}, that is, if $X\in\g$ is $\ad$-semisimple in
$\g$, then it is $\ad$-semisimple in the Lie algebra of $\GGGG$.
The details are similar to the proof of Lemma \ref{L:linear} and
are omitted here.
\end{remark}


\section{Proof of Theorem \ref{T:main1}}

Based on the preparations of the previous two sections, in this
section we prove the closedness of $\Ad$-semisimple conjugacy
classes in connected Lie groups and $\ad$-semisimple adjoint
orbits in real Lie algebras.

\begin{theorem}\label{T:closed}
Let $G$ be a connected Lie group with Lie algebra $\g$. Then\\
(i) $\Ad$-semisimple conjugacy classes in $G$ are closed embedded
submanifold of $G$;\\
(ii) $\ad$-semisimple adjoint orbits in $\g$ are closed embedded
submanifold of $\g$.
\end{theorem}

\begin{proof}
Let $G'=G/Z(G)_0$, where $Z(G)_0$ is the identity component of the
center $Z(G)$ of $G$. Let $\pi:G\rightarrow G'$ be the quotient
homomorphism. Then the adjoint representation $\Ad:G\rightarrow
GL(\g)$ induces naturally a representation $\rho:G'\rightarrow
GL(\g)$, such that $\rho\circ\pi=\Ad$. Note that $\ker(\rho)$ is
discrete in $G'$.

(i) Let $C$ be an $\Ad$-semisimple conjugacy class in $G$. Then
$C'=\pi(C)$ is a conjugacy class in $G'$. Since all elements of
$\rho(C')=\Ad(C)$ are semisimple in $GL(\g)$, by Corollary
\ref{C:polynomial}, $C'$ is a closed embedded submanifold of $G'$.
So $$M=\pi^{-1}(C')=C\cdot Z(G)_0$$ is a closed embedded
submanifold of $G$.

Now consider the transitive action of $G\times Z(G)_0$ on the
manifold $M$, defined by $$(h,z).x=hxh^{-1}z.$$ Choose a $g\in
C\subset M$, and let $L\subset G\times Z(G)_0$ be the isotropic
group of $g$. Then the map $$\varphi:(G\times Z(G)_0)/L\rightarrow
M$$ defined by $$\varphi((h,z)L)=hgh^{-1}z$$ is a diffeomorphism.

By Theorem \ref{T:finite}, $\Gamma(C)=g^{-1}C\cap Z(G)$ is a
finite subgroup of $Z(G)$. So $\Gamma_0(C)=g^{-1}C\cap
Z(G)_0=\Gamma(C)\cap Z(G)_0$ is a finite subgroup of $Z(G)_0$.
Then $Z(G)_0/\Gamma_0(C)$ is a Lie group. Let $$\alpha:G\times
Z(G)_0\rightarrow Z(G)_0/\Gamma_0(C)$$ be the epimorphism defined
by $$\alpha(h,z)=[z],$$ where $[z]$ is the image of $z$ under the
quotient homomorphism $Z(G)_0\rightarrow Z(G)_0/\Gamma_0(C)$. For
$(h,z)\in L$, $hgh^{-1}z=(h,z).g=g$, so
$z^{-1}=g^{-1}hgh^{-1}\in\Gamma_0(C)$, and then $[z]$ is trivial
in $Z(G)_0/\Gamma_0(C)$. This shows that $L\subset\ker(\alpha)$.
Then $\alpha$ induces a smooth map $$\widetilde{\alpha}:(G\times
Z(G)_0)/L\rightarrow Z(G)_0/\Gamma_0(C)$$ defined by
$$\widetilde{\alpha}((h,z)L)=[z].$$ It is obvious that $(G\times
Z(G)_0)/L$ is a fiber bundle with base space $Z(G)_0/\Gamma_0(C)$,
fiber type $\widetilde{\alpha}^{-1}([e])$, and projection
$\widetilde{\alpha}$.

We claim that $$\widetilde{\alpha}^{-1}([e])=\varphi^{-1}(C).$$
Firstly, let $(h,z)L\in\widetilde{\alpha}^{-1}([e])$. Then
$[z]=[e]$, that is, $z\in\Gamma_0(C)$. So there exists $k\in G$
such that $z=g^{-1}kgk^{-1}$. Then
$$\varphi((h,z)L)=hgh^{-1}z=hgzh^{-1}=hg(g^{-1}kgk^{-1})h^{-1}=(hk)g(hk)^{-1}\in
C,$$ that is, $(h,z)L\in\varphi^{-1}(C)$. Conversely, let
$(h',z')L\in\varphi^{-1}(C)$. Then there exists $k'\in G$ such
that $\varphi((h',z')L)=h'gz'h'^{-1}=k'gk'^{-1}$. This implies
$z'=g^{-1}(h'^{-1}k')g(h'^{-1}k')^{-1}$. So $z'\in\Gamma_0(C)$.
Hence $(h',z')L\in\widetilde{\alpha}^{-1}([e])$. This verifies the
claim.

As the fiber above $[e]$,
$\varphi^{-1}(C)=\widetilde{\alpha}^{-1}([e])$ is a closed
embedded submanifold of $(G\times Z(G)_0)/L$. Since $\varphi$ is a
diffeomorphism, $C$ is a closed embedded submanifold of $M$, hence
a closed embedded submanifold of $G$. Item (i) is proved.

(ii) Let $O$ be an $\ad$-semisimple adjoint orbit in $\g$. Then
$O'=d\pi(O)$ is an adjoint orbit in $\g'$, the Lie algebra of
$G'$. Since all elements of $d\rho(O')=\ad(O)$ is semisimple in
$\gl(\g)$, by Corollary \ref{C:polynomial}, $O'$ is a closed
embedded submanifold of $\g'$. So $$N=(d\pi)^{-1}(O')=O+Z(\g)$$ is
a closed embedded submanifold of $\g$.

Consider the transitive action of $G\times Z(\g)$ on the manifold
$N$, defined by $$(h,Y).W=\Ad(h)W+Y.$$ Choose an $X\in O\subset
N$, and let $K\subset G\times Z(\g)$ be the isotropic group of
$X$. Then the map $$\psi:(G\times Z(\g))/K\rightarrow N$$ defined
by $$\psi((h,Y)K)=\Ad(h)X+Y$$ is a diffeomorphism.

Let $$\beta:G\times Z(\g)\rightarrow Z(\g)$$ be the projection to
the second factor. For $(h,Y)\in K$, $\Ad(h)X+Y=(h,Y).X=X$, so
$-Y=-X+\Ad(h)X\in\gamma(O)$. But by Theorem \ref{T:finite},
$\gamma(O)$ is trivial. So $Y=0$. This shows that
$K\subset\ker(\beta)$. Then $\beta$ induces a smooth map
$$\widetilde{\beta}:(G\times Z(\g))/K\rightarrow Z(\g)$$ defined by
$$\widetilde{\beta}((h,Y)K)=Y.$$ It is obvious that $(G\times
Z(\g))/K$ is a fiber bundle with base space $Z(\g)$, fiber type
$\widetilde{\beta}^{-1}(0)$, and projection $\widetilde{\beta}$.
Similar to the proof of (i), we have
$\widetilde{\beta}^{-1}(0)=\psi^{-1}(O)$.

As the fiber above $0\in Z(\g)$,
$\psi^{-1}(O)=\widetilde{\beta}^{-1}(0)$ is a closed embedded
submanifold of $(G\times Z(\g))/K$. Since $\psi$ is a
diffeomorphism, $O$ is a closed embedded submanifold of $N$, hence
a closed embedded submanifold of $\g$. This proves (ii).
\end{proof}


\section{Proof of Theorem \ref{T:main2}}

We give the proof of Theorem \ref{T:main2} in this section. We
first prove a lemma.

\begin{lemma}\label{L:covering}
Let $\pi:G\rightarrow G'$ be a covering homomorphism of connected
Lie groups. If $C'$ is an $\Ad$-semisimple conjugacy class in
$G'$, then every connected component of $\pi^{-1}(C')$ is a
conjugacy class in $G$.
\end{lemma}

\begin{proof}
By Theorem \ref{T:closed}, $C'$ is a closed embedded submanifold
of $G'$. Let $\widetilde{C}$ be a connected component of
$\pi^{-1}(C')$. Then $\widetilde{C}$ is a closed embedded
submanifold of $G$. Since $G$ is connected, $\widetilde{C}$ is
invariant under the conjugation of $G$. Let $C$ be a conjugacy
class in $G$ contained in $\widetilde{C}$. By Theorem
\ref{T:closed}, $C$ is a closed embedded submanifold of $G$, hence
a closed embedded submanifold of $\widetilde{C}$. But $\dim C=\dim
C'=\dim\widetilde{C}$. By the connectedness of $\widetilde{C}$, we
must have $C=\widetilde{C}$.
\end{proof}

\begin{remark}
Lemma \ref{L:covering} does not hold without the
$\Ad$-semisimplicity assumption.
\end{remark}

\begin{theorem}\label{T:intersection}
Let $\alpha:H\rightarrow G$ be a homomorphism of connected Lie
groups. Suppose $\ker(\alpha)$ is discrete. Let the Lie algebras
of $G$ and $H$
be $\g$ and $\ah$, respectively. We have\\
(i) If $C$ is an $\Ad$-semisimple conjugacy class in $G$, then
every connected component of
$\alpha^{-1}(C)$ is a conjugacy class in $H$;\\
(ii) If $O$ is an $\ad$-semisimple adjoint orbit in $\g$, then
every connected component of $(d\alpha)^{-1}(O)$ is an adjoint
orbit in $\ah$.
\end{theorem}

\begin{proof}
(i) We first observed that $\alpha^{-1}(C)$ is invariant under the
conjugation of $H$. So $\alpha^{-1}(C)$ is the union of a family
of conjugacy classes in $H$. But the connectedness of $H$ implies
that conjugacy classes in $H$ are connected. So every connected
component of $\alpha^{-1}(C)$ is the union of a family of
conjugacy classes in $H$. We prove that every connected component
of $\alpha^{-1}(C)$ is a conjugacy class in $H$. The proof is
divided into three steps.

\emph{Step (a)}. We prove (i) under the additional assumptions
that $\alpha$ is injective and $\Gamma(C)$ is trivial. In this
case, $H$ can be identified with $\alpha(H)$, which is a Lie
subgroup of $G$, and then $\alpha^{-1}(C)$ is identified with
$C\cap H$. Note that under this identification, the prior topology
on $H$ may be different from the subspace topology on $H$ induced
from $G$. We call the prior topology on $H$ the
\emph{$H$-topology}, and call a connected component of $C\cap H$
with respect to the $H$-topology an \emph{$H$-connected component}
of $C\cap H$.

Let $C_i$ be an $H$-connected component of $C\cap H$. Consider the
adjoint homomorphism $\Ad_G=\Ad:G\rightarrow \Ad(G)$. Then
$\Ad_G(C_i)\subset\Ad_G(C)\cap\Ad_G(H)$. Let $C'_i$ be the
$\Ad_G(H)$-connected component of $\Ad_G(C)\cap\Ad_G(H)$
containing $\Ad_G(C_i)$. Since all elements of $\Ad_G(C)$ are
semisimple, by Corollary \ref{C:polynomial}, the conjugacy class
$\Ad_G(C)$ in $\Ad(G)$ is an $\Ad(G)$-connected component of
$Z=\{A\in\Ad(G)|f(A)=0\}$, where $f$ is the minimal polynomial of
$\Ad_G(h_0)$ for some $h_0\in H$. So $C'_i$ is an
$\Ad_G(H)$-connected component of $Z\cap\Ad_G(H)$. By Corollary
\ref{C:polynomial} again, we conclude that $C'_i$ is a conjugacy
class in $\Ad_G(H)$. Let $g_1,g_2\in C_i$. Then $\Ad_G(g_1),
\Ad_G(g_2)\in C'_i$, and then there exists $h\in H$ such that
$\Ad_G(g_2)=\Ad_G(h)\Ad_G(g_1)\Ad_G(h)^{-1}$. So $g_2=hg_1h^{-1}z$
for some $z\in Z(G)$. But $g_1$ and $g_2$ are conjugate in $G$. So
there is $g\in G$ such that $g_2=gg_1g^{-1}$. This implies
$gg_1g^{-1}=hg_1h^{-1}z=hg_1zh^{-1}$. Hence
$z=g_1^{-1}(h^{-1}g)g_1(h^{-1}g)^{-1}\in g_1^{-1}C\cap
Z(G)=\Gamma(C)$. But we have assumed that $\Gamma(C)$ is trivial.
So $z=e$, and then $g_2=hg_1h^{-1}$. This shows that $C_i$ is a
conjugacy class in $H$.

\emph{Step (b)}. We prove (i) under the additional assumption that
$\alpha$ is injective. As we have done in step (a), we identify
$H$ with $\alpha(H)$. Let $G'=G/\Gamma(C)$. By Theorem
\ref{T:finite}, $\Gamma(C)$ is finite. So the quotient
homomorphism $\pi:G\rightarrow G'$ is a covering homomorphism. In
particular, $\pi|_H:H\rightarrow\alpha(H)$ is a covering
homomorphism. Let $C_i$ be an $H$-connected component of $C\cap
H$, and let $C'$ be a conjugacy class in $H$ contained in $C_i$.
Then $\pi(C')$ is a conjugacy class in $\pi(H)$, and we have
$\pi(C')\subset\pi(C_i)\subset\pi(C)\cap\pi(H)$. Let $C'_i$ be the
$\pi(H)$-connected component of $\pi(C)\cap\pi(H)$ containing
$\pi(C_i)$. By Lemma \ref{L:invariant},
$\Gamma_{G'}(\pi(C))=\pi(\Gamma(C))$ is trivial. So by step (a),
$C'_i$ is a conjugacy class in $\pi(H)$ containing $\pi(C')$. This
forces $\pi(C')=\pi(C_i)=C'_i$. Hence $C'\subset
C_i\subset(\pi|_H)^{-1}(C'_i)$. By Lemma \ref{L:covering}, $C'$ is
an $H$-connected component of $(\pi|_H)^{-1}(C'_i)$. As an
$H$-connected subset of $(\pi|_H)^{-1}(C'_i)$ containing $C'$,
$C_i$ must coincide with $C'$. So $C_i$ is a conjugacy class in
$H$.

\emph{Step (c)}. We finish the proof of item (i). Let $C_i$ be a
connected component of $\alpha^{-1}(C)$, and let $C'_i$ be the
$\alpha(H)$-connected component of $C\cap\alpha(H)$ containing
$\alpha(C_i)$. Then $C_i$ is a connected component of
$\alpha^{-1}(C'_i)$. But by step (b), $C'_i$ is a conjugacy class
in $\alpha(H)$. So by Lemma \ref{L:covering}, $C_i$ is a conjugacy
class in $H$.

(ii) Since $d\alpha$ is injective, $\ah$ can be viewed as a
subalgebra of $\g$. We want to prove that if $O\cap\ah$ is
nonempty, then every connected component of $O\cap\ah$ is an
adjoint orbit in $\ah$. Let $O_i$ be a connected component of
$O\cap\ah$. Then $\ad_\g(O_i)\subset\ad_\g(O)\cap\ad_\g(\ah)$. Let
$O'_i$ be the connected component of $\ad_\g(O)\cap\ad_\g(\ah)$
containing $\ad_\g(O_i)$. Since all elements of $\ad_\g(O)$ are
semisimple, by Corollary \ref{C:polynomial}, the adjoint orbit
$\ad_\g(O)$ in $\ad(\g)$ is a connected component of
$\az=\{B\in\ad(\g)|p(B)=0\}$, where $p$ is the minimal polynomial
of $\ad_\g(Y_0)$ for some $Y_0\in\ah$. So $O'_i$ is a connected
component of $\az\cap\ad_\g(\ah)$. By Corollary \ref{C:polynomial}
again, we conclude that $O'_i$ is an adjoint orbit in
$\ad_\g(\ah)$. Let $X_1,X_2\in O_i$. Then $\ad_\g(X_1),
\ad_\g(X_2)\in O'_i$, and then there exists $h\in H$ such that
$\ad_\g(X_2)=\Ad(\Ad_G(h))\ad_\g(X_1)$. So $X_2=\Ad_G(h)X_1+Y$ for
some $Y\in Z(\g)$. But $X_1$ and $X_2$ lie in the same adjoint
orbit in $\g$. So there is $g\in G$ such that $X_2=\Ad_G(g)X_1$.
This implies $Y=-\Ad_G(h)X_1+\Ad_G(g)X_1\in\gamma(O)$. By Theorem
\ref{T:finite}, $\gamma(O)$ is trivial. So $Y=0$, and then
$X_2=\Ad_G(h)X_1+Y$. This shows that $O_i$ is an adjoint orbit in
$\ah$. (ii) is proved.
\end{proof}

\end{document}